\documentclass[11pt]{article}

\usepackage{amsmath, amsfonts, amssymb}
\usepackage{lipsum}
\usepackage{graphicx}
\usepackage{epstopdf}
\usepackage{bm}
\usepackage{xcolor}
\usepackage{multirow} 
\usepackage{url}
\usepackage{hyperref} 

\newcommand{\bd}[1]{\mathbf{#1}}

\DeclareMathOperator\erf{erf}
\DeclareMathOperator\erfc{erfc}

\title{Localized evaluation and fast summation in the extrapolated regularization method for integrals in Stokes flow}
\author{Joseph Siebor\thanks{Department of Mathematics, Farmingdale State College, SUNY, Farmingdale, NY 11735, USA siebjp@farmingdale.edu} 
\and Svetlana Tlupova\thanks{Department of Mathematics, Farmingdale State College, SUNY, Farmingdale, NY 11735, USA tlupovs@farmingdale.edu}}

\date{\today}

\begin{document}

\maketitle

\begin{abstract} 
Boundary integral equation methods are widely used in the solution of many partial differential equations. The kernels that appear in these surface integrals are nearly singular when evaluated near the boundary, and straightforward numerical integration produces inaccurate results. In Beale and Tlupova (Adv. Comput. Math, 2024), an extrapolated regularization method was proposed to accurately evaluate the nearly singular single and double-layer surface integrals for harmonic potentials or Stokes flow. The kernels are regularized using a smoothing parameter, and then a standard quadrature is applied. The integrals are computed for three choices of the smoothing parameter to find the extrapolated value to fifth order accuracy. In this work, we apply several techniques to reduce the computational cost of the extrapolated regularization method applied to the Stokes single and double layer integrals. First, we use a straightforward OpenMP parallelization over the target points. Second, we note that the effect of the regularization is local and evaluate only the local component of the sum for three values of the smoothing parameter. The non-local component of the sum is only evaluated once and reused in the other sums. This component is still the computational bottleneck as it is $O(N^2)$, where $N$ is the system size. We apply the kernel-independent treecode to these far-field interactions to reduce the CPU time. We carry out experiments to determine optimal parameters both in terms of accuracy and efficiency of the computations. We then use these techniques to compute Stokes flow around two spheres that are nearly touching.
\end{abstract}



\section{Introduction}

Boundary integral equation methods are widely used in the solution of many partial differential equations. For example, the Stokes single and double layer surface integrals on a closed surface $\Gamma$ are,
\begin{equation}
\label{Stokes_SL}
	u_i(\bd{y}) = \frac{1}{8\pi}\int_\Gamma S_{ij}(\bd{y,x}) f_j(\bd{x}) dS(\bd{x}),
\end{equation}
\begin{equation}
\label{Stokes_DL}
	v_i(\bd{y}) = \frac{1}{8\pi}\int_\Gamma T_{ijk} (\bd{y,x}) q_j(\bd{x}) n_k(\bd{x})dS(\bd{x}),
\end{equation}
where $f$ and $q$ are density functions, and
\begin{equation}
\label{Stokeslet}
	S_{ij}(\bd{y,x}) = \frac{\delta_{ij}}{|\bd{y} - \bd{x}|} + \frac{(y_i - x_i)(y_j - x_j)}{|\bd{y} - \bd{x}|^3},
\end{equation}
\begin{equation}
\label{Stresslet}
	T_{ijk} (\bd{y,x}) = -\frac{6(y_i - x_i)(y_j - x_j)(y_k - x_k)}{|\bd{y} - \bd{x}|^5},
\end{equation}
are the Stokeslet and the stresslet kernels, $r = |\bd{y} - \bd{x}|$, $\delta_{ij}$ is the Kronecker delta, and $i,j,k = 1,2,3$. 

The kernels that appear in surface integrals such as~\eqref{Stokeslet},\eqref{Stresslet} are singular when evaluated on the boundary and nearly singular when evaluated near the boundary. Straightforward numerical integration produces inaccurate results. The accurate evaluation of nearly singular integrals has been a very active area of research for the past couple of decades, and several approaches are available in the literature. Among them are an analytical resolution of
the singularity by a coordinate change~\cite{bruno-kunyansky-01}, interpolation procedures using values at points away from the singularity~\cite{ying-biros-zorin-06, zorincplx}, quadrature by expansion~\cite{baggetorn, klinttorn,
klockner-barnett-greengard-oneil-13, siegel-tornberg-18}, singularity subtraction~\cite{helsing-13}, a corrected trapezoidal rule~\cite{nitsche}, asymptotic expansions of the kernel~\cite{carvalho18}, expansions of the density rather than the kernel to reduce the singularities~\cite{perez}, and regularization with corrections~\cite{b04, byw, tlupova-beale-19}. 

In~\cite{extrap}, an extrapolated regularization method was proposed to accurately evaluate the nearly singular single and double layer surface integrals, for harmonic potentials or Stokes flow. The kernels of the form $1/r^p$, where $r=|\hat{\bd{y}}|$ and $p=1,3,5$, are replaced with regularized versions $s_i(r/\delta)/r^p$, $i=1,2,3$, with $\delta$ being the regularization, or smoothing, parameter. Using local analysis, an equation is derived for the error due to regularization. To eliminate the leading terms in this error, the integrals are computed with three choices of the smoothing parameter. Then a system of three equations is solved for an extrapolated value with the error reduced to $O(\delta^5)$ uniformly on or near the surface. The method can be extended to $O(\delta^7)$ accuracy by evaluating the regularized integrals for four choices of $\delta$ and solving a system of four equations to find the extrapolated value. A standard quadrature is applied to evaluate the integrals, with the discretization error negligible compared to the regularization error if $\delta/h$ is kept large enough. For constant $\delta/h$ and moderate resolution $h$ on the surface, a combined error of about $O(h^5)$ is observed for the case of three $\delta$'s.  The extrapolated regularization approach of~\cite{extrap} is simple and direct in the sense that the work required is similar to that for a surface integral with a smooth integrand, except that three (or four) related integrals must be computed rather than one. No special gridding or separate treatment of the singularity is needed. Geometric information about the surface is not needed other than normal vectors. 

The cost of evaluating the discretized sums is $O(MN)$ where $N$ is the number of quadrature (or source) points and $M$ is the number of evaluation (or target) points. Evaluating these sums for three choices of the regularization parameter $\delta$ will then add a factor of three to the CPU time. The goal of this paper is to demonstrate several approaches that can be stacked to greatly reduce the cost of the method. First, we note that the effect of regularization is local, that is, kernels such as $S$ and $T$ in~\eqref{Stokeslet},\eqref{Stresslet} are nearly indistinguishable from their regularized versions away from the singularity. We then propose to evaluate only the local components of the sum for the three $\delta$ values, and not use regularization in the well-separated components. This removes the additional "extrapolation" factor in the CPU estimate. However, $O(MN)$ is still too expensive for large systems. For efficiency, fast summation methods suitable for regularized kernels~\cite{rbfstokes, wang-krasny-tlupova-20, yingradial} could be used, and we apply the kernel-independent treecode~\cite{wang-krasny-tlupova-20} in this work to reduce the cost to $O(M\log N)$. Finally, we use OpenMP to parallelize the outer loop over the target points and to reduce the overall CPU time to about $O(\frac{1}{n}M\log N)$, where $n$ is the number of processors. 

We summarize the extrapolated regularization method of~\cite{extrap} in Section~\ref{section:extrap_reg}. Our strategies for reducing the overall computational cost of the method are described in Section~\ref{section:local_tree}. Numerical examples that demonstrate the effectiveness of these approaches are presented in Section~\ref{section:results}.


\section{Extrapolated regularization}
\label{section:extrap_reg}

First, subtraction is used in the integrals~\eqref{Stokes_SL} and~\eqref{Stokes_DL},
\begin{equation}
\label{Stokes_SL1}
	u_i(\bd{y}) = \frac{1}{8\pi}\int_\Gamma S_{ij}(\bd{y,x})\, [f_j(\bd{x}) - f_k(\bd{x}_0)n_k(\bd{x}_0)n_j(\bd{x})] \, dS(\bd{x}),
\end{equation}
\begin{equation}
\label{Stokes_DL1}
	v_i(\bd{y}) = \frac{1}{8\pi}\int_\Gamma T_{ijk} (\bd{y,x})\, [q_j(\bd{x}) - q_j(\bd{x}_0)] \, n_k(\bd{x})\, dS(\bd{x}) + \chi(\bd{y})\, q_i(\bd{x}_0),
\end{equation}
where $\bd{x}_0$ is the closest point on $\Gamma$, and $\chi(\bd{y})=1,0,1/2$ when $\bd{y}$ is inside, outside, or on $\Gamma$, respectively. In addition, as explained in~\cite{extrap}, the stresslet kernel is rewritten as follows,
\begin{equation}
\label{Tsplit}
	T_{ijk} =  T_{ijk}^{(1)} + T_{ijk}^{(2)} = - 6\left(\frac{t_{ijk}^{(1)}}{r^3} + \frac{t_{ijk}^{(2)} - (r^2 - b^2)\, t_{ijk}^{(1)}}{r^5}\right),
\end{equation}
with
\begin{equation}
	t_{ijk}^{(1)}  = bn_in_jn_k - (\hat{x}_i n_j n_k + n_i\hat{x}_j n_k + n_i n_j \hat{x}_k),
\end{equation}
\begin{equation}
	t_{ijk}^{(2)} = b(\hat{x}_i\hat{x}_j n_k + \hat{x}_i n_j\hat{x}_k + n_i\hat{x}_j\hat{x}_k) - \hat{x}_i\hat{x}_j\hat{x}_k,
\end{equation}
where $n_i$ and $\hat{x}_i$ are the $i$th components of $\bd{n}_0$ and $\hat{\bd{x}} = \bd{x}-\bd{x}_0$.

We compute~\eqref{Stokes_SL1} and~\eqref{Stokes_DL1} as 
\begin{equation}
\label{Stokes_SL2}
	u^\delta_i(\bd{y}) = \frac{1}{8\pi}\int_\Gamma S^\delta_{ij}(\bd{y,x})\, [f_j(\bd{x}) - f_k(\bd{x}_0)n_k(\bd{x}_0)n_j(\bd{x})] \, dS(\bd{x}),
\end{equation}
\begin{equation}
\label{Stokes_DL2}
	v^\delta_i(\bd{y}) = \frac{1}{8\pi}\int_\Gamma T^\delta_{ijk} (\bd{y,x}) \, [q_j(\bd{x}) - q_j(\bd{x}_0)] \, n_k(\bd{x})\, dS(\bd{x}) + \chi(\bd{y})\, q_i(\bd{x}_0),
\end{equation}
with $S_{ij}$ and $T_{ijk}$ replaced with the regularized versions,
\begin{equation}
\label{Sreg} 
	S_{ij}^\delta(\bd{y,x}) = \frac{\delta_{ij}}{r}\, s_1(r/\delta) + \frac{(y_i - x_i)(y_j - x_j)}{r^3}\, s_2(r/\delta)\,,
\end{equation}
\begin{equation}
\label{Treg} 
	T_{ijk}^\delta = T_{ijk}^{(1)}\, s_2(r/\delta) + T_{ijk}^{(2)}\, s_3(r/\delta),
\end{equation}
where 
\begin{equation}
\label{s1}
    s_1(r) = \erf(r) = \frac{2}{\sqrt{\pi}}\int_0^r e^{-s^2}\,ds,
\end{equation}
\begin{equation}
\label{s2}
    s_2(r) = \erf(r) - \frac{2}{\sqrt{\pi}}\, r \, e^{-r^2},
\end{equation}
\begin{equation}
\label{s3}
    s_3(r) = \erf(r) - \frac{2}{\sqrt{\pi}}\left(r + \frac23 r^3 \right) e^{-r^2}.
\end{equation}
Since $\erf(r/\delta) \to 1$ rapidly as $r/\delta$ increases, this regularization has a localized effect. This is discussed in more detail in Sec.~\ref{section:local_tree}.

The regularization error, namely, $u^\delta - u$ or $v^\delta-v$, is typically $O(\delta)$.
If $\bd{y}$ is near the surface, then we write $\bd{y} = \bd{x}_0 + b\bd{n}$, where $\bd{x}_0$ is the closest point on $\Gamma$, $\bd{n}$ is the outward normal vector at $\bd{x}_0$, and $b$ is the signed distance.  From a series expansion for $\bd{x}$ near $\bd{x}_0$ and $b$ near $0$, it was shown in~\cite{extrap} that
\begin{equation}
	\label{extrap_eqn}
	u_i(\bd{y}) + c_1 \rho I_0(\lambda) + c_2\rho^3 I_2(\lambda) = u_i^\delta (\bd{y}) + O(\delta^5),
\end{equation}
uniformly for $\bd{y}$ near the surface,
where $\lambda = b/\delta$, $\rho = \delta/h$, $c_1$ and $c_2$ are unknown coefficients, and $I_0$ and $I_2$ are integrals occurring in the derivation that are known explicitly,
\begin{equation}
	\label{I0}
	I_0(\lambda) = e^{-\lambda^2}/\sqrt{\pi} - |\lambda|\erfc|\lambda|,
\end{equation}
\begin{equation}
	\label{I2}
	I_2(\lambda) = \frac23\left((\frac12 - \lambda^2) e^{-\lambda^2}/\sqrt{\pi} + |\lambda|^3 \erfc|\lambda| \right),
\end{equation}
where $\erfc = 1 - \erf$. It is important that $c_1, c_2$ do not depend on $\delta$ or $\lambda$.
To obtain an accurate value of $u_i$, we calculate the regularized integrals $u_i^\delta$ for three different choices of $\delta$, with the same grid size $h$, resulting in a system of three equations with three unknowns.  We can then solve for the exact integral $u_i$ within error $\delta^5$. We typically choose $\delta_i = \rho_i h$ with $\rho_i = 3, 4, 5$. This procedure, and the equations~\eqref{extrap_eqn}-\eqref{I2}, are identical for the stresslet integral $v_i^\delta$ in~\eqref{Stokes_DL2}.

\subsection{Evaluation on the surface}
\label{subsec:eval_surf}

When the integrals need to be evaluated on the surface, e.g., when solving integral equations, the extrapolated regularization can be used as described, but a more direct method would be to replace the smoothing factors $s$ in~\eqref{s1}-\eqref{s3} with modified functions~\cite{tlupova-beale-19, novel},
\begin{equation}
    s_1^\sharp(r) = \erf(r) + \frac{2}{3\sqrt{\pi}} \left(5r - 2r^3\right) e^{-r^2},
\end{equation}
\begin{equation} 
    s_2^\sharp(r) = \erf(r) - \frac{2}{3\sqrt{\pi}} \left(3r - 14r^3 + 4r^5\right) e^{-r^2},
\end{equation}
\begin{equation}
    s_3^\sharp(r) = \erf(r) - \frac{2}{9\sqrt{\pi}} \left(9r + 6r^3 -4r^5\right) e^{-r^2}.
\end{equation}
Using these functions to regularize the kernels will result in $O(\delta^5)$ accuracy directly without the need for extrapolation.

\subsection{Quadrature}
\label{subsec:quadrature}

We briefly discuss the quadrature rule for surface integrals introduced in~\cite{wilson-10, byw} and used in~\cite{extrap}. The surface is covered with a three-dimensional grid with spacing $h$. A set $\Gamma_3$ of quadrature points have the form $\bd{x} = (ih,jh,x_3)$ such that their projections on the $(x_1,x_2)$ plane are grid points, and $|\bd{n}(\bd{x}) \cdot \bd{e}_3| \geq \cos\theta$, where $\bd{n}(\bd{x})$ is the unit normal at $\bd{x}$, and we take $\theta = 70^o$. Sets $\Gamma_1$ and $\Gamma_2$ are defined similarly. For a function $f$ on the surface the integral is computed as
\begin{equation} 
    \label{quadrature}
    \int_\Gamma f(\bd{x})\,dS(\bd{x}) \,\approx\, \sum_{i=1}^3 \sum_{\bd{x} \in \Gamma_i} f(\bd{x})w_i(\bd{x}) \,,
    \qquad w_i(\bd{x}) = \frac{\psi_i(\bd{n}(\bd{x}))}{|n_i(\bd{x})|} \,h^2. 
\end{equation}
The weights $w_i$ are determined by a partition of unity $\psi_1, \psi_2, \psi_3$ on the unit sphere that is applied to the normal vector $\bd{n} = (n_1,n_2,n_3)$ at each point.  This is done by using the bump function
\begin{equation} 
    b(r) = \exp(2 r^2/(r^2-1))\,, \quad |r|<1\,; \qquad
    b(r) =  0\,, \quad |r| \geq 1,              
\end{equation} 
and defining
\begin{equation} 
    \beta_i(\bd{n}) = b(\cos^{-1}|n_i|)/\theta) \,, \quad \psi_i(\bd{n}) = \beta_i(\bd{n})/\left(\sum_{j=1}^3 \beta_j(\bd{n})\right).  
\end{equation}
The quadrature rule~\eqref{quadrature} has high order accuracy as allowed by the smoothness of the surface and the integrand. The weights cut off the sum in each plane, and each sum has the character of the trapezoidal rule without boundary; see~\cite{wilson-10}. For the full discussion of the error we refer to~\cite{extrap, neglect, byw}. We only note here that the quadrature error is small for $\delta/h\geq 2$ and decreases rapidly as $\delta/h$ increases. For moderate resolution sizes, setting $\delta=\rho h$, with $\rho=2,3,4$ or $\rho=3,4,5$, should work well in most cases, resulting in about $O(h^5)$ convergence. To ensure that the regularization error dominates the discretization error for small $h$ we can choose $\delta$ proportional to $h^q$, with $q < 1$, so that $\delta/h$ increases as $h \to 0$.


\section{Localized evaluation and fast summation using treecode}
\label{section:local_tree}

\subsection{OpenMP parallelization}

To reduce the computational cost of evaluating the discrete sums at $M$ target points, we use OpenMP to split the calculations over $n$ threads. The threads parallelize the calculations and have access to a shared memory space. The speedup gained via OpenMP is theoretically a factor of $n$ threads used.

\subsection{Localizing regularization to the near field}
\label{subsec:loc_reg}

\begin{figure}\centering
\includegraphics[scale=0.45]{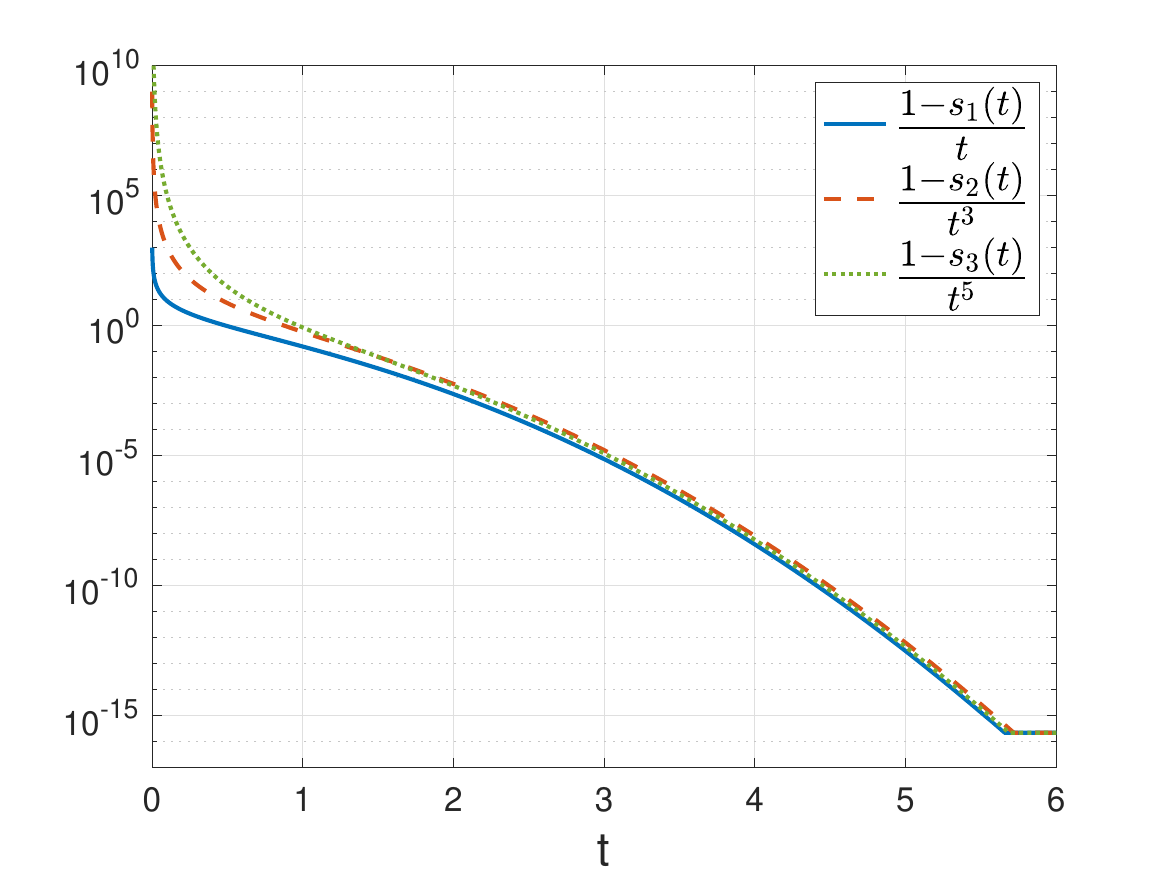}
\caption{The localized effect of smoothing using $s_1$, $s_2$, $s_3$.}
\label{fig:s1s2}
\end{figure}

Next, we focus on reducing the computational cost of evaluating the discrete sums at a single target point, which is $O(N)$ where $N$ is the number of quadrature points. As described in Sec.~\ref{section:extrap_reg}, extrapolated regularization relies on evaluating the surface integrals for three choices of the regularization parameter $\delta_1, \delta_2, \delta_3$ in order to obtain an accurate extrapolated value. The numerical tests in~\cite{extrap} showed that setting $\delta_i=\rho_i h$, $\rho_i=3,4,5$, seems to be a reliable choice for the Stokes integrals in most scenarios. This adds a factor of three to the computational complexity, which we symbolically denote as $3\cdot O(N)$ for each target point. However, the effect of regularization is localized, with the difference between the singular and regularized versions of the kernels decreasing very rapidly. As seen in Fig.~\ref{fig:s1s2}, the values of $(1-s_1(t))/t$, $(1-s_2(t))/t^3$, and $(1-s_3(t))/t^5$ diminish to machine precision for $t>6$, so we will simply set $s_1=s_2=s_3=1$ for $t=r/\delta > 6$. 

\subsection{Localizing extrapolated regularization}
\label{subsec:loc_extrap}

In light of the above observation, we denote the discretized integrals  by
\begin{equation}
    \label{near_far}
    u_{\delta} = u_{\delta}^{near} + u^{far}, \\[6pt]
\end{equation}
where $u^{near}$ and $u^{far}$ denote the components of the sum with quadrature points that are near and far from the target point, respectively. For the near component, regularization is used as described above for $r/\delta \leq 6$, and for the far component, regularization is not necessary. We then localize the extrapolated regularization by reusing the well-separated component of the sum $u^{far}$ for all three values of $\delta$. This is a very simple modification that effectively removes the factor of three in the CPU time. The computational complexity is now $O(N+2N_{near})$, where $O(2N_{near})$ accounts for the overhead of evaluating the near components. The smaller the "near radius", the smaller the overhead. At the same time, the near radius must be wide enough not to affect the accuracy.

\subsection{Fast summation using the kernel-independent treecode}
\label{subsec:treecode}

\begin{figure}[htb]\centering
\includegraphics[scale=0.6]{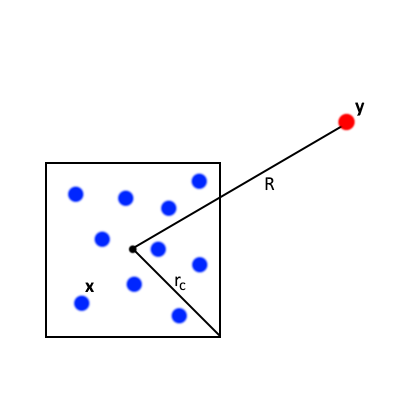}
\caption{A target point $\bd{y}$ and a source cluster $C$. $R$ denotes the target-cluster distance, and $r_c$ is the cluster radius.}
\label{target-cluster}
\end{figure}

The computational bottleneck in evaluating the discretized integrals comes from the far-field interactions, that is, the computation of $u^{far}$ in~\eqref{near_far}. A suitable fast summation method, e.g.,~\cite{rbfstokes, wang-krasny-tlupova-20, yingradial}, can be used to speed up this component. We use the kernel-independent treecode~\cite{wang-krasny-tlupova-20}, which is based on barycentric Lagrange interpolation at Chebyshev points to approximate the well-separated particle-cluster interactions. This treecode requires only kernel evaluations and is expected to reduce the $O(N)$ sum to $O(\log N)$ for each target point. In order to efficiently determine which quadrature points are near and which are not, the entire set of quadrature points is first partitioned recursively into a hierarchical octree structure, dividing each cluster into eight child clusters until no more than a user-defined number $N_0$ of points remains in each cluster, called a leaf. The algorithm then cycles recursively through clusters, rather than individual quadrature points. For a given target point $\bd{y}$ and a cluster of radius $r_c$ (see Fig.~\ref{target-cluster}), we define a target-cluster radius $R$ (distance from $\bd{y}$ to the center of the cluster). The target and the cluster are well-separated if a multipole acceptance criterion (MAC) is met,
\begin{equation}
    \label{MAC}
    \frac{r_c}{R} \leq \theta,
\end{equation}
where $\theta$ is a user-defined parameter.

Suppose for a target point $\bd{y}$ and a well-separated cluster of source particles $C=\{\bd{x}_j\}$ the interaction is denoted by
\begin{equation}
    \label{particle-cluster}
    u(\bd{y},C) = \sum_{\bd{x}_j\in C} K(\bd{y},\bd{x}_j)\, g(\bd{x}_j),
\end{equation}
where $K$ refers to parts of the Stokeslet $S_{ij}$ or the stresslet $T_{ijk}$, and $g$ is the corresponding weight function (given by the density functions $f$ or $q$, and the quadrature weights $w$). The cluster $C$ is a rectangular box whose sides are aligned with the coordinate axes, and Chebyshev points $\bd{s}_\bd{k}$, $\bd{k}=(k_1,k_2,k_3)$ with $k_l=0,...,p$, are mapped to this box in each of the $x_1,x_2,x_3$ directions. Using a tensor product in three directions, the interaction~\eqref{particle-cluster} is computed using the kernel approximation
\begin{equation}
    K(\bd{y},\bd{x}) \approx \sum_{k_1=0}^p \sum_{k_2=0}^p \sum_{k_3=0}^p K(\bd{y},\bd{s}_\bd{k}) L_{k_1}(x_1) L_{k_2}(x_2) L_{k_3}(x_3), 
\end{equation}
where $L_l$ are polynomials of degree $p$ in each variable $x_1,x_2,x_3$ obtained using the barycentric Lagrange form. Using this approximation, the well-separated particle-cluster interaction~\eqref{particle-cluster} is computed as
\begin{equation}
    \label{particle-cluster-KITC}
    u(\bd{y},C) \approx \sum_{k_1=0}^p \sum_{k_2=0}^p \sum_{k_3=0}^p K(\bd{y},\bd{s}_\bd{k})\, \hat{g_\bd{k}},
\end{equation}
where
\begin{equation}
    \label{mod_weights}
    \hat{g_\bd{k}} = \sum_{\bd{x}_j\in C} L_{k_1}(x_{j1}) L_{k_2}(x_{j2}) L_{k_3}(x_{j3}) \, g(\bd{x}_j)
\end{equation}
are referred to as modified weights. The speedup in the computations is due to the fact that the modified weights in~\eqref{mod_weights} for a given cluster are independent of the target $\bd{y}$, so they can be precomputed and reused for different targets. Only kernel evaluations $K(\bd{y},\bd{s}_\bd{k})$ at the Chebyshev points $\bd{s}_\bd{k}$ are needed in~\eqref{particle-cluster-KITC}, making the algorithm kernel-independent. 

To apply the method to evaluating the far-field sums in the Stokeslet~\eqref{Stokes_SL1}, we pre-compute the modified weights as in~\eqref{mod_weights} for $f_j(\bd{x}) w(\bd{x})$ and $n_j(\bd{x})w(\bd{x})$, $j=1,2,3$ (6 modified weights in total). Similarly, for the stresslet~\eqref{Stokes_DL1}, we compute the modified weights for $q_j(\bd{x})n_k(\bd{x}) w(\bd{x})$ and $n_k(\bd{x}) w(\bd{x})$, $j,k=1,2,3$ (12 modified weights in total). The modified weights computed this way require only a small fraction of the total CPU time. We refer to~\cite{wang-krasny-tlupova-20} for details in computing the modified weights and the estimates for time and storage. Since we are evaluating the well-separated component of the sums, the kernels $S_{ij}$ and $T_{ijk}$ are used in their original forms~\eqref{Stokeslet} and~\eqref{Stresslet} without regularization. 

The overall algorithm then proceeds as follows. The user-specified parameters are $N_0$, the leaf size, $\theta$, the MAC parameter, and $p$, the degree of the Lagrange polynomial approximation. First, the quadrature points are divided into the tree of clusters. Then, the modified weights are computed and stored. We loop over the target points using OpenMP. For each target point, the algorithm cycles through the clusters, starting with the root cluster and proceeding to the child clusters. For a given target point and source cluster, we check the MAC~\eqref{MAC}. If it is satisfied, then the target and cluster are well-separated and we compute the interaction using the approximation in~\eqref{particle-cluster-KITC} which uses the original kernels without regularization~\eqref{Stokeslet} and~\eqref{Stresslet}. If the MAC is not satisfied, the code checks each child of the cluster. If the cluster is a leaf, the target-cluster interaction is computed by a direct sum using the regularized kernels~\eqref{Sreg} and~\eqref{Treg} for three values of the regularization parameter $\delta$, using the strategy of Sec.~\ref{subsec:loc_reg}.


\section{Numerical results}
\label{section:results}

We begin our numerical tests by calculating the Stokeslet integral near a translating spheroid. We use OpenMP to parallelize the outer loop over the target points and localize the regularization as in Sec.~\ref{subsec:loc_reg}. Our tests verify that doing so preserves the original accuracy while giving the expected speedup of almost $n$ in the CPU time, where $n$ is the number of processes. Next, we localize the extrapolated regularization by reusing the well-separated component of the sums in the computation of all three $u_\delta$ as in Sec.~\ref{subsec:loc_extrap}. We use the kernel-independent treecode of~\cite{wang-krasny-tlupova-20} to reduce the well-separated component of the sum from $O(N)$ to $O(\log N)$. For $M$ target points, the overall estimate is then $O(\frac{1}{n}M\log N)$. Finally, we use these techniques to compute Stokes flow around two spheres that are nearly touching. All tests were performed on a MacBook Pro Sequoia 15.5 OS with a 2.3 GHz Quad-Core Intel Core i7 processor. The code was written in double precision C\texttt{++} and compiled using the Clang compiler with the -O2 optimization.


\subsection{Parallel evaluation using OpenMP}

We present tests of the CPU time by calculating the Stokes velocity near a spheroid $x_1^2 + 4x_2^2 + 4x_3^2 = 1$ translating with velocity $(1,0,0)$. The flow velocity at any point $\bd{y}$ outside the spheroid is determined by the single layer integral~\eqref{Stokes_SL} with 
\begin{equation}
	f_1(\bd{x}) = \frac{F_0}{\sqrt{1-3x_1^2/4}},\quad f_2(\bd{x})=f_3(\bd{x}) = 0,
\end{equation}
and $F_0$ is a constant, see~\cite{extrap}. The exact solution is given in~\cite{chwang, liron, tlupova-beale-19}; it has a maximum amplitude and an $L^2$ norm of about 1. We consider a regular three-dimensional grid with size $h$, and evaluate the Stokeslet integral at grid points outside the spheroid that are a distance of at most $h$ to the surface. We use extrapolated regularization with three values of $\delta=\rho h$, $\rho=3,4,5$, leading to about $O(h^5)$ accuracy for moderate resolution. The errors reported are the maximum errors and the $L^2$ errors, defined as
\begin{equation}
	\label{error_L2}
	||e||_2 = \left( \sum_{\bd{y}} |e(\bd{y})|^2 / N \right)^{1/2},
\end{equation}
where $e(\bd{y})$ is the error at $\bd{y}$ and $N$ is the number of evaluation points. The errors were reported in Fig. 10 of~\cite{extrap}. 

Table~\ref{OpenMP} shows the number of quadrature and target points, the errors and CPU times for different grid sizes $h$. The sums were evaluated directly in both serial and parallel mode with 4 threads. In the sum for each $\delta$, the regularized kernels were used only for $r/\delta \leq 6$, as discussed in Sec.~\ref{subsec:loc_reg}. The errors we obtain in Table~\ref{OpenMP} match those in~\cite{extrap}, with OpenMP speeding up the computations by a factor of about 3.6. We use OpenMP in all remaining tests.


\begin{table}[!htb]
\centering
\begin{tabular}{|c|c|c|c|c||c|c|c|}
\hline
\multirow{2}{*}{$1/h$} & \multirow{2}{*}{\# quads} & \multirow{2}{*}{\# targets} & \multirow{2}{*}{max err} & \multirow{2}{*}{$L^2$ err} & \multicolumn{3}{c|}{CPU time (sec.)} \\ 
 & & & & & serial & parallel, 4 threads & speedup \\ 
\hline
32 & 6958 & 5856 & 3.27e-3 & 3.35e-4 & 3.56 & 1 & 3.56 \\
\hline
64 & 27934 & 22720 & 2.03e-4 & 1.65e-5 & 39 & 11 & 3.55 \\
\hline
128 & 112006 & 89684 & 8.09e-6 & 6.02e-7 & 566 & 155 & 3.65 \\
\hline
256 & 448094 & 356032 & 2.71e-7 & 1.95e-8 & 8797 & 2461 & 3.57 \\
\hline
\end{tabular}
\caption{Stokeslet integral on a prolate spheroid, at grid points within distance $h$ outside the spheroid. CPU times with serial and parallel computations.}
\label{OpenMP}
\end{table}


\subsection{Fast summation using treecode}

Here we present the treecode performance with the error values in Table~\ref{OpenMP} as a benchmark, and propose a strategy for choosing the optimal treecode parameters. There are three parameters we must set. First, $\theta$ is a measure of a target point and a cluster being well-separated. A larger value of $\theta$ will lead to more treecode approximations, thereby increasing the error and decreasing the CPU time. Then, $N_0$ is the maximum number of points in any leaf in the tree. The effect of $N_0$ on the CPU time is less straightforward. Smaller $N_0$ means there are more levels in the tree, and recursive cycling through these levels takes more time, but with smaller clusters, the well-separation criterion will be satisfied more often which would save CPU also. Finally, $p$ is the degree of the polynomial used in the treecode approximations. Increasing $p$ will lead to higher accuracy but is slower, as the CPU estimate includes a factor of $p^3$. As we refine the grid size $h$, the number of quadrature, or source, points will increase within a fixed geometric space, while the error is expected to decrease as $O(h^5)$, and our goal is to find values for these parameters that will speed up the computations while preserving the accuracy with respect to decreasing $h$.

\begin{figure}[htb]
\centering
\includegraphics[scale=0.275]{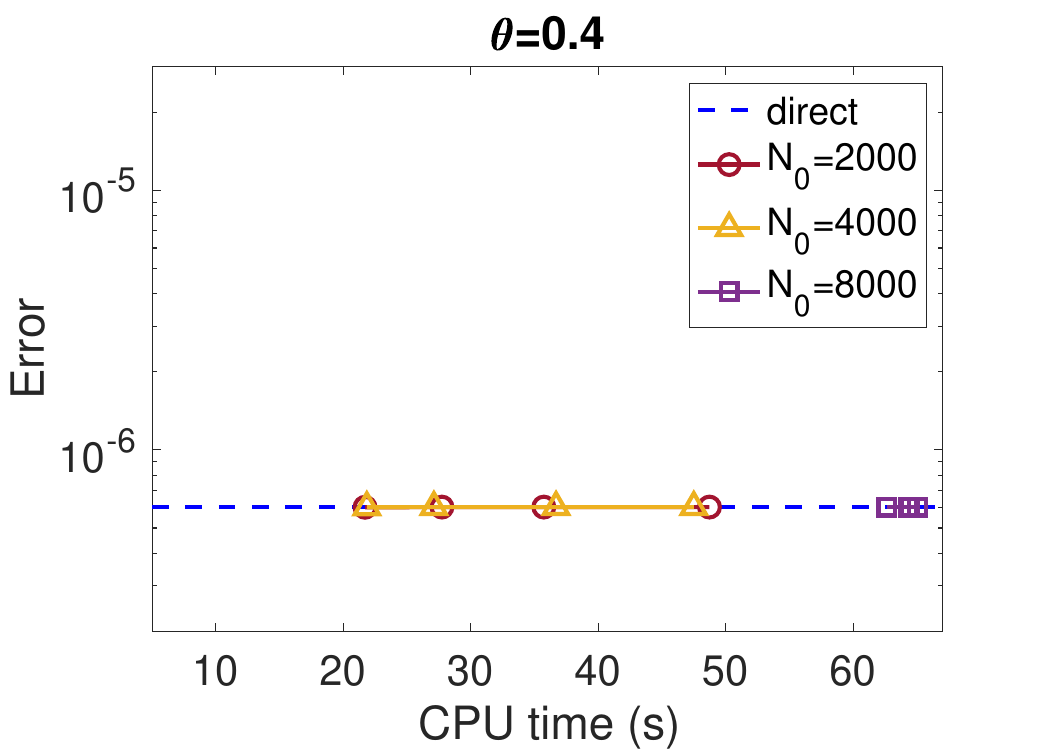}
\includegraphics[scale=0.275]{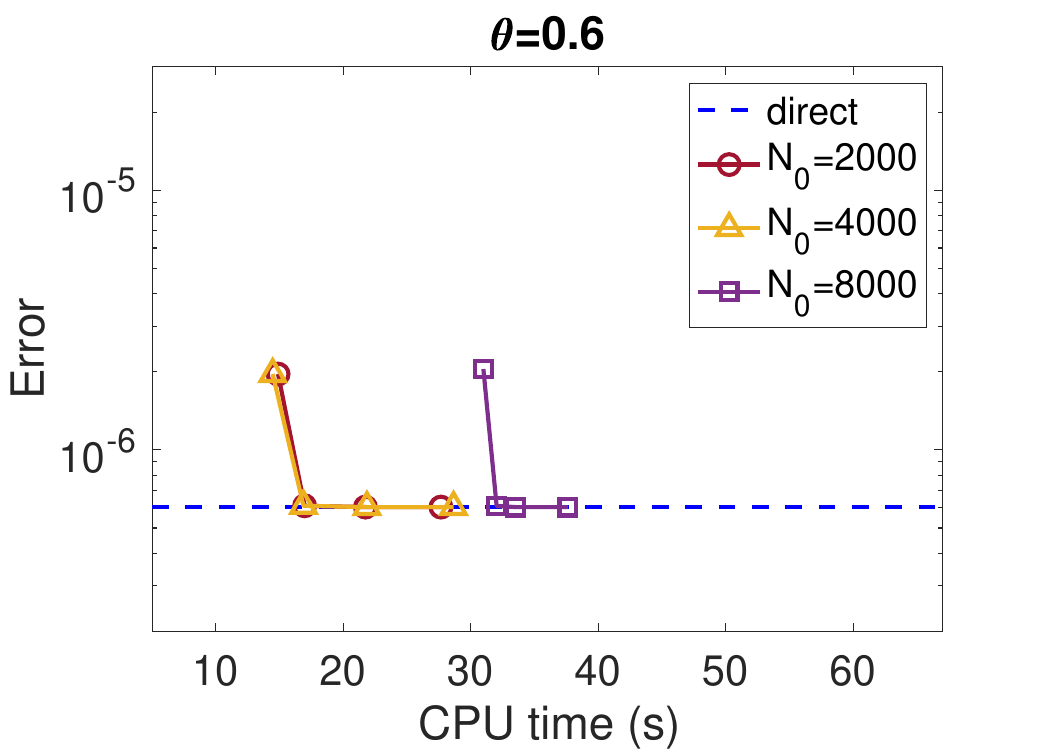}
\includegraphics[scale=0.275]{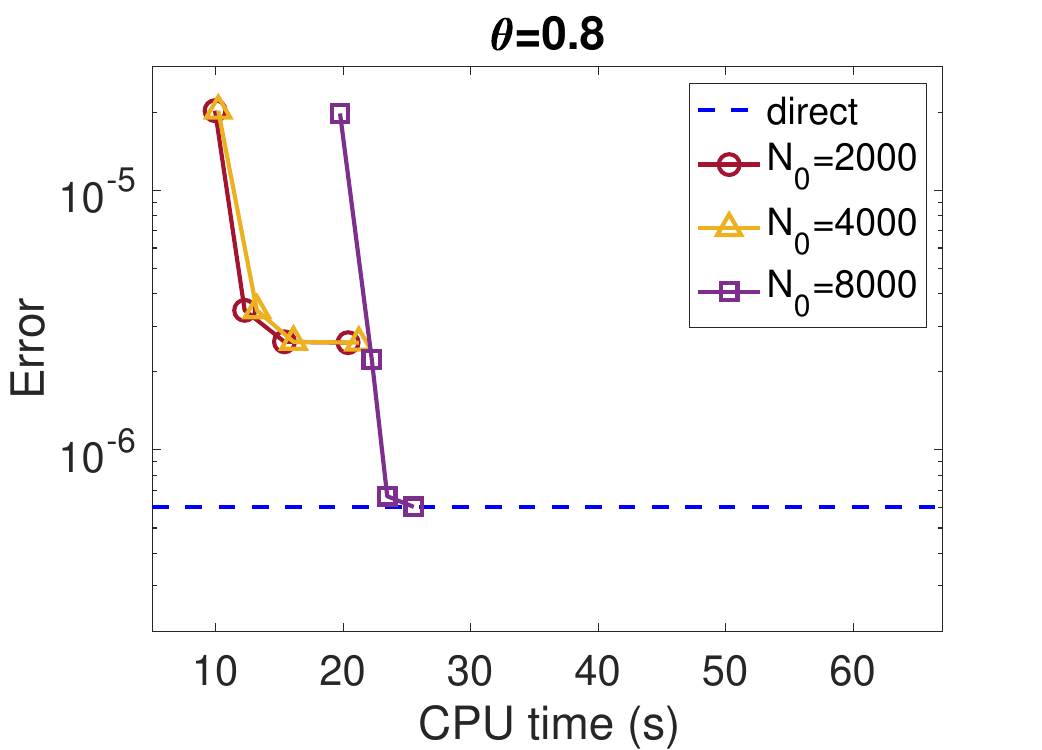}
\caption{Stokeslet integral at grid points within $h$ of the translating spheroid, $h=\frac{1}{128}$; $L^2$ error vs. treecode CPU time for $\theta=0.4$ (left), $\theta=0.6$ (middle), $\theta=0.8$ (right); degree $p=6,8,10,12$ (increasing from left to right). The direct sum is 155 sec (blue horizontal line).}
\label{theta_choice}
\end{figure}

We start our experiments with $\theta$. We fix $h=1/128$, and compute the solution outside a translating spheroid as above using the treecode approximations on the original Stokeslet kernel $S_{ij}$ when the target and cluster are well-separated, that is, when the MAC~\eqref{MAC} is met for a given $\theta$. This is the far component of the sums $u^{far}$, and is reused for all three values of $\delta$, as described in Sec.~\ref{subsec:loc_extrap}. When the MAC is not met, each child of the cluster is checked, and when the cluster is a leaf, the near field direct summation is done using the regularized kernel $S^\delta_{ij}$ (again, regularization applied for $r/\delta \leq 6$ only). Figure~\ref{theta_choice} shows the error and CPU time for three values of $\theta=0.4, 0.6, 0.8$. For each $\theta$, we used different values of $N_0=\{2K, 4K, 8K\}$ and $p=\{6,8,10,12\}$. We compare the error to the benchmark that is the direct sum. With the small value $\theta=0.4$, all computations were accurate but slow. With the large value $\theta=0.8$, most computations were not as accurate as the direct sum. For a moderate value of $\theta=0.6$ however, all computations for $p=8$ or above met the accuracy threshold while being more efficient than the small $\theta$ results. For this reason, we fix $\theta=0.6$ as our optimal value in all the remaining tests.

\begin{figure}[htb]
\centering
\includegraphics[scale=0.375]{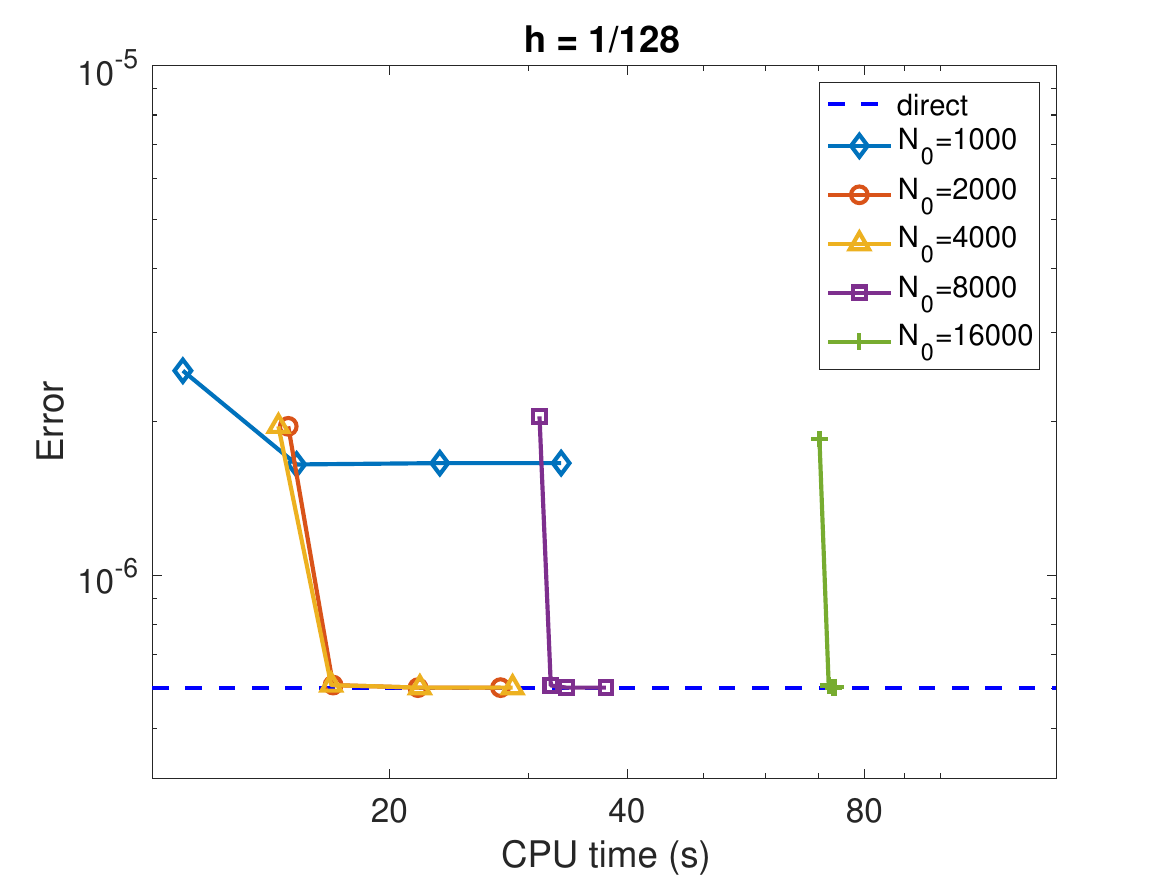}
\includegraphics[scale=0.375]{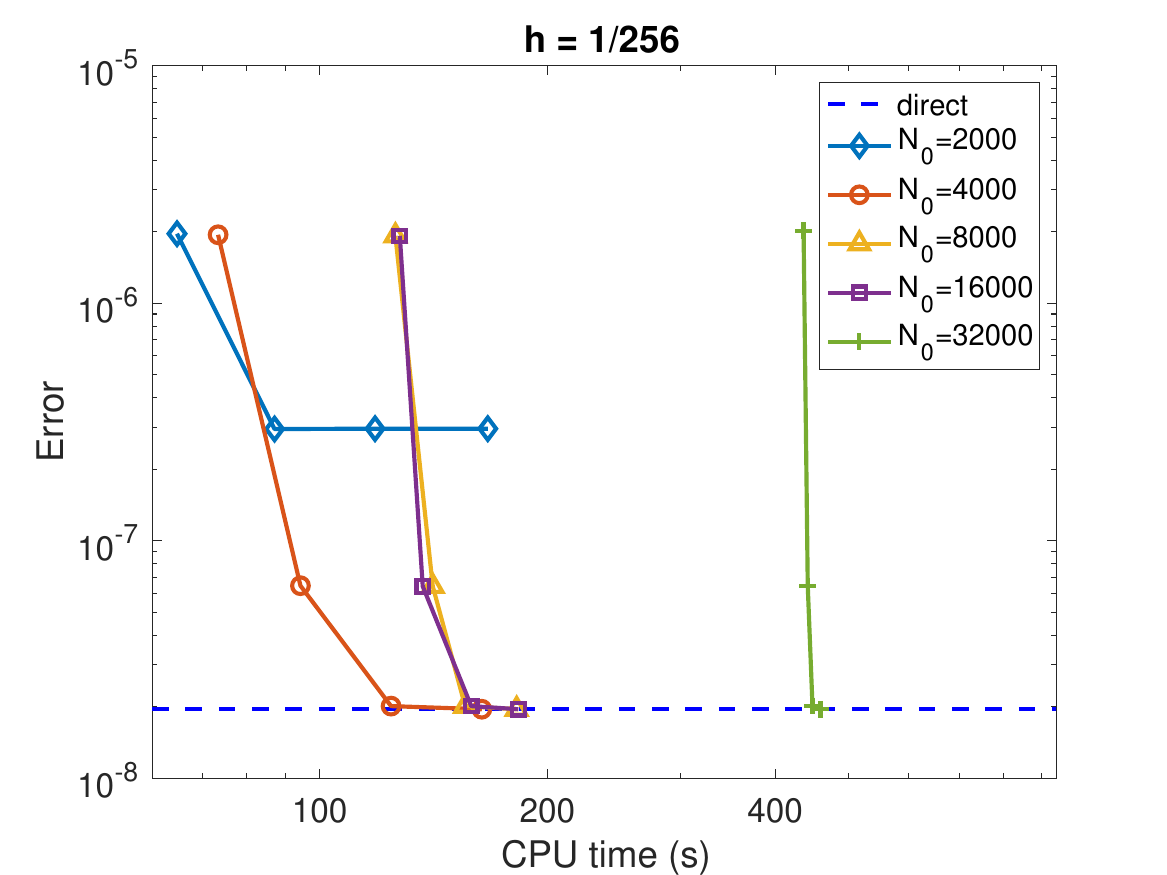}
\caption{Stokeslet integral at grid points within $h$ of the translating spheroid, $h=\frac{1}{128}$ (left), $h=\frac{1}{256}$ (right); $L^2$ error vs. treecode CPU time for $\theta=0.6$, degree $p=6,8,10,12$ (increasing from left to right). The direct sums are 155 sec and 2461 sec (blue horizontal lines).}
\label{N0_choice}
\end{figure}

With $\theta=0.6$ fixed, we then analyze the effect of $N_0$ and $p$ on the error and CPU time while $h$ decreases from $h=1/128$ to $h=1/256$, shown in Fig.~\ref{N0_choice}. For both values of $h$, we see that having too small $N_0$ (in blue) will lead to larger error that will not decrease with $p$. At the same time, having too large $N_0$ (in green) will be inefficient. These two "floor" and "ceiling" values of $N_0$ seem to double as $h$ decreases by a factor of 2. For $h=1/128$, the values are $N_0=1K$ and $N_0=16K$, while for $h=1/256$, the values are $N_0=2K$ and $N_0=32K$. We then choose $N_0=4K$ as the optimal value in the first case, and double it to $N_0=8K$ for the smaller $h$. We will continue to double $N_0$ as $h$ decreases by a factor of 2. Regarding the value of $p$, the error threshold is met at $p=8$ and $p=10$ for the two $h$ values, and we propose that $p$ is increased by an additional 2 each time $h$ is decreased by a factor of 2. To summarize, our proposed strategy is:
\begin{equation}
\label{optimal_params}
	\theta=0.6;\qquad h = \frac{1}{64}: \ \ N_0=2000, \ p=6; \qquad h\to h/2: \ \ N_0\to 2N_0, \ p\to p+2.
\end{equation}
For $h$ less than $1/64$, we will fix $N_0=1K$ and $p=6$ to ensure small error, as the CPU time is small in those cases.

\begin{figure}[htb]
\centering
\includegraphics[scale=0.5]{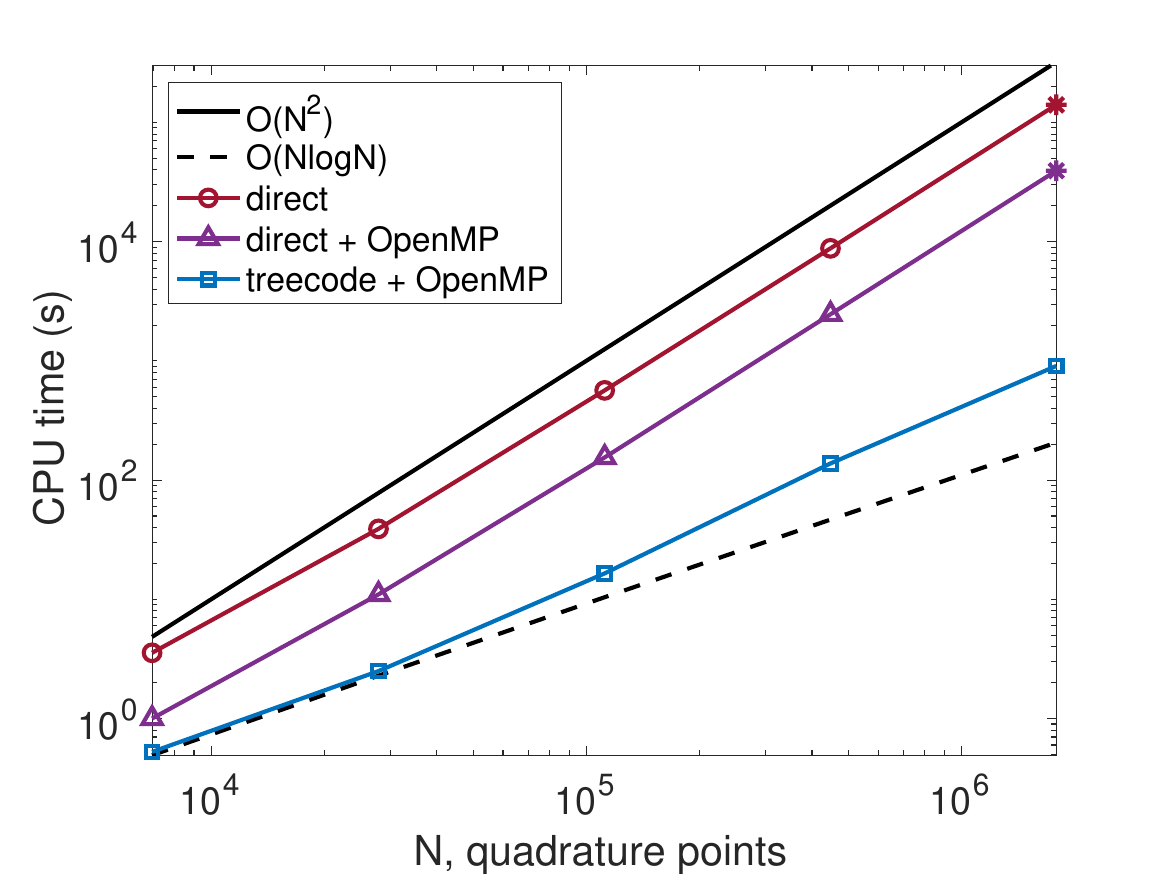}
\caption{Stokeslet integral at grid points within $h$ of the translating spheroid, CPU time vs. $N$, the number of quadrature points (estimated time is marked by *). The optimal parameters according to~\eqref{optimal_params} were used in the treecode.}
\label{results}
\end{figure}

Finally, Fig.~\ref{results} shows the CPU time as the system size $N$ grows with grid refinement. In the treecode, optimal parameters as set by equation~\ref{optimal_params} were used for each $N$. The direct sum is clearly $O(N^2)$, while the treecode is about $O(N\log N)$, with a slight increase due to varying $N_0$ and increasing $p$ (as $p$ is the degree of the interpolating polynomial in each dimension, the CPU estimate includes a factor of $p^3$).

\subsection{Two drops nearly touching in Stokes flow}

Here we consider an example of Stokes flow around two spheres that present an interface of fluids of different viscosities. The interface velocity satisfies an integral equation~\cite{pozbook},
\begin{align}
	\label{IE-2Surf}
	(\lambda_b+1) u_i(\textbf{x}_0) = &-\frac{1}{4\pi \mu_0}\sum_{m=1}^{2}\int_{\partial\Omega_m} S_{ij}(\textbf{x}_0,\textbf{x}) [f]_j(\textbf{x})dS(\textbf{x}) \nonumber\\
	&+ \sum_{m=1}^{2}\frac{\lambda_m-1}{4\pi}\int_{\partial\Omega_m} T_{ijk} (\textbf{x}_0,\textbf{x}) u_j(\textbf{x}) n_k(\textbf{x})dS(\textbf{x})
\end{align}
for $\textbf{}{x}_0\in \partial\Omega_b$, $b=1,2$, $\lambda_b=\mu_b/\mu_0$. We take the external viscosity $\mu_0=1$ and internal viscosities $\mu_1 = \mu_2 = 2$, so that the viscosity ratios for each interface are $\lambda_1=\lambda_2=2$. In~\ref{IE-2Surf}, $[f]$ is the discontinuity in the surface force given by $[\bd{f}] = 2\gamma H \bd{n} - \nabla_S\gamma$, where $\gamma$ is the surface tension, $H$ is the mean curvature, and $\bd{n}$ is the outward unit normal. We define the surface tension $\gamma = 1+(x_1-x_c)^2$, where $x_c$ is the $x$-coordinate of the center of the sphere. We center one sphere at the origin and the other at $(2,0,\epsilon)$, where $\epsilon = 1/16^3$. 

For each $b=1,2$, there are four integrals to evaluate, two Stokeslet integrals and two stresslet integrals. For each of these, one is an integral over the self-surface, where $\bd{x}_0\in \partial\Omega_b$ and we integrate over $\Omega_b$, and the other is the integral over the nearby surface. The former involves on-surface integration and we use the special regularization described in Sec.~\ref{subsec:eval_surf} with $\delta=3h$. The latter involves nearly singular integration and extrapolated regularization is used with $\delta=\{3h, 4h, 5h\}$. We apply OpenMP, localized regularization, and the treecode to all the integrals. As before, for the nearly singular integrals, the localized regularization is evaluated through direct summation for three values of the smoothing parameter $\delta$ so that extrapolation can be done. The optimal treecode parameters outlined in~\eqref{optimal_params} were used for each $h$.

We solve the integral equation~\ref{IE-2Surf} using GMRES with a tolerance $10^{-10}$. Since the exact solution is not known, we check the convergence rates by defining
\begin{equation}
	\label{error_h}
	e_h (\bd{x}) = \bd{u}_{h} (\bd{x}) - \bd{u}_{h/2} (\bd{x}),
\end{equation}
and taking either the max or the $L^2$ norm of this error as defined by~\eqref{error_L2}, where $N$ is the number of surface points given by $h$, the larger of the two grid sizes used. 

\begin{figure}[htb]
\centering
\includegraphics[scale=0.6]{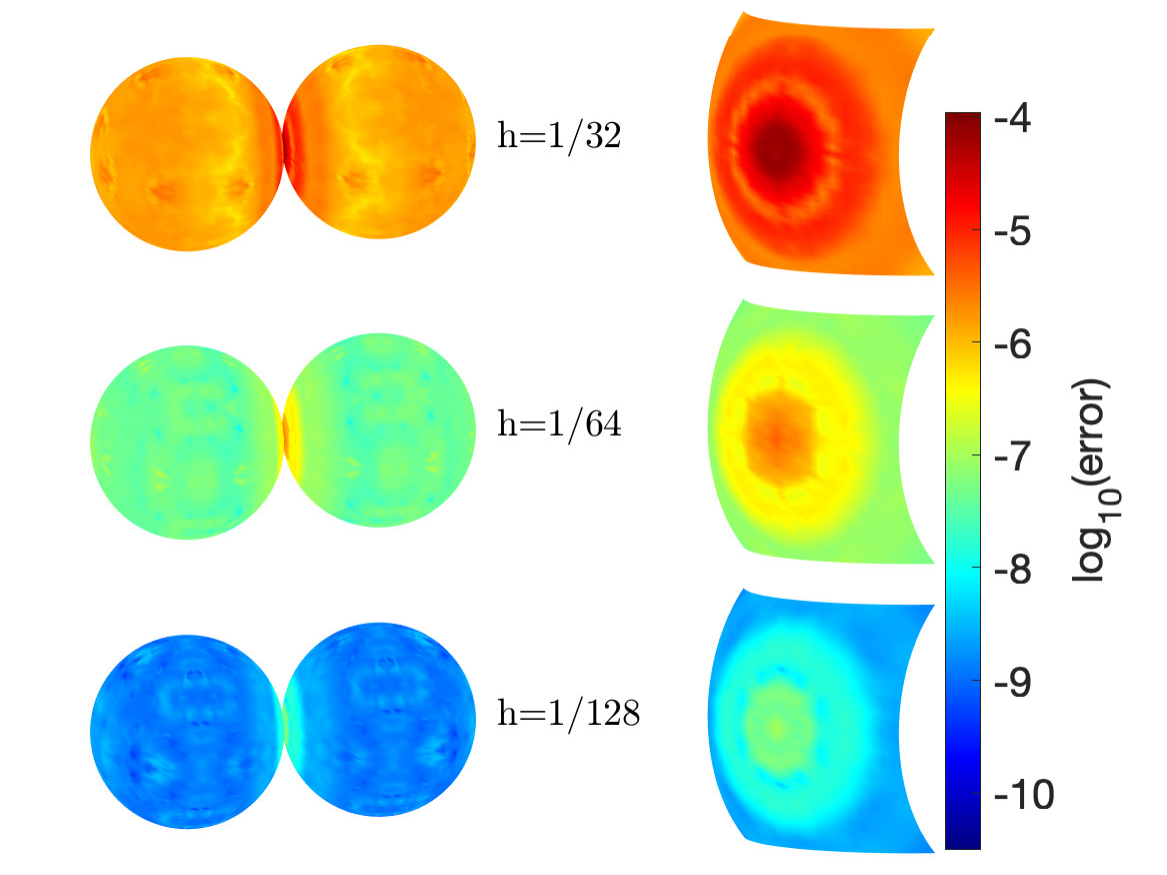}
\caption{Two drops nearly touching in Stokes flow. Error (log scale) in the solution of the integral equation~\eqref{IE-2Surf} on the surface of the drops, with the close-ups of the near singular region (right).}
\label{fig:2Surf}
\end{figure}

Table~\ref{table:2Surf} shows the number of GMRES iterations it took to reach the tolerance, the solution error and the convergence rate. For each $h$, it took 12 iterations to converge, and the order of convergence is about $O(h^5)$, verifying that the choice of the treecode parameters has not affected the expected high order of convergence of the extrapolated regularization method. Figure~\ref{fig:2Surf} shows the error on the surface of both drops for different grid sizes $h$, with the close ups of the region where the drops are nearly touching. While this region has the largest errors on the surface, the error decreases quickly with grid refinement.
\begin{table}[htb]
\centering
\begin{tabular}{|c|c|c|c|c|c|c|c|}
\hline
$1/h$ & \ $N_0$ \ & \ $p$ \ & \# of iterations & max err & order & $L^2$ err & order \\ 
\hline
16 & 1000 & 6 & 13 & 3.17e-3 & -- & 1.19e-4 & -- \\
\hline
32 & 1000 & 6 & 12 & 1.07e-4 & 4.89 & 6.11e-6 & 4.28 \\
\hline
64 & 2000 & 8 & 12 & 4.06e-6 & 4.72 & 2.11e-7 & 4.86 \\
\hline
128 & 4000 & 10 & 12 & 1.02e-7 & 5.31 & 6.37e-9 & 5.05 \\
\hline
256 & 8000 & 12 & 12 & -- & -- & -- & -- \\
\hline
\end{tabular}
\caption{Two drops nearly touching in Stokes flow. Solution of the integral equation~\eqref{IE-2Surf} on the surface of the drops. Treecode parameters $\theta=0.6$, $N_0$ and $p$; the number of GMRES with tolerance $10^{-10}$; the max and $L^2$ errors in the solution and order of convergence.}
\label{table:2Surf}
\end{table}

\section*{Acknowledgments}
This work was supported in part by the National Science Foundation grant DMS-2012371.

\bibliographystyle{plain}
\bibliography{Bib}

\end{document}